\documentclass{article}
\usepackage[utf8]{inputenc}
\usepackage{amsthm}
\usepackage{amsmath}
\usepackage{amsfonts}
\usepackage{babel}
\usepackage{geometry}
\usepackage{fullpage}
\usepackage[colorlinks=false,hidelinks]{hyperref}
\usepackage{xurl}
\usepackage{color}
\usepackage{graphicx}
\usepackage[dvipsnames]{xcolor}
\usepackage[shortlabels]{enumitem}
\usepackage[nosort,capitalise]{cleveref}
\crefname{enumi}{}{}
\crefname{equation}{}{}
\usepackage{tikz}
\usetikzlibrary{matrix}
\usepackage{booktabs}

\definecolor{AfonsoBlue}{RGB}{30,65,123}

\title{Randomstrasse101: Open Problems of 2025}

\author{
  Afonso S. Bandeira\thanks{ASB: Department of Mathematics, ETH Zürich, Rämistrasse 101, 8092 Zurich, Switzerland. \texttt{bandeira@math.ethz.ch}} 
  \and
  Daniil Dmitriev\thanks{DD: Department of Statistics and Data Science at the Wharton School, University of Pennsylvania, PA, USA. \texttt{dmitrievdaniil97@gmail.com}}
  \and
  Kevin Lucca\thanks{KL: Department of Mathematics, ETH Zürich, Rämistrasse 101, 8092 Zurich, Switzerland. \texttt{kevin.lucca@ifor.math.ethz.ch}}
  \and
  Petar Nizi\'{c}-Nikolac\thanks{PNN: Department of Mathematics, ETH Zürich, Rämistrasse 101, 8092 Zurich, Switzerland. \texttt{petar.nizic-nikolac@ifor.math.ethz.ch}}
  \and
  Almut Rödder\thanks{AR: Department of Mathematics, ETH Zürich, Rämistrasse 101, 8092 Zurich, Switzerland. \texttt{almutmagdalena.roedder@ifor.math.ethz.ch}}
}

\date{\today}

\usepackage{natbib}
\usepackage{graphicx}

\newcommand{\EE}{\mathbb{E}}
\newcommand{\RR}{\mathbb{R}}
\newcommand{\CC}{\mathbb{C}}
\newcounter{sectionforenv}

\newtheorem{theorem}{Theorem}[sectionforenv]

\makeatletter
\@addtoreset{sectionforenv}{section}
\makeatother
\newtheorem{problem}{Problem}
\newtheorem{definition}[theorem]{Definition}

\newtheorem{conjecture}[problem]{Conjecture}
\newtheorem{openproblem}[problem]{Open Problem}

\usepackage{tocloft}

\cftsetindents{section}{1em}{5em}
\begin{document}

\addtocounter{section}{7}
\addtocounter{problem}{15}

\maketitle

\begin{abstract}
\texttt{Randomstrasse101} is a blog dedicated to Open Problems in Mathematics, with a focus on Probability Theory, Computation, Combinatorics, Statistics, and related topics. This manuscript serves as a stable record of the Open Problems posted in 2025, with the goal of easing academic referencing. The blog can currently be accessed at~\texttt{randomstrasse101.math.ethz.ch}.
\end{abstract}


\section*{Introduction}

In this manuscript we include the blog entries in
\texttt{Randomstrasse101} of 2025, containing a total of sixteen Open Problems. The 2024 version can be found here~\cite{RandomstrasseProblems2024}. \texttt{Randomstrasse101} is a blog created and maintained by our group at the Department of Mathematics at ETH Z\"{u}rich\footnote{The inspiration for the name of the blog will be clear to the reader after looking up the department's address and recalling the Probability Theory focus of the blog.}.

The focus is on mathematical open problems in Probability Theory, Computation, Combinatorics, Statistics, and related topics. The goal is not to necessarily write about the most important open questions in these fields but simply to discuss open questions and conjectures\footnote{Conjectures here should be interpreted as mathematical statements that we do not know not to be true, and for which a proof or a refutal would be interesting progress. We do not necessarily have very high confidence that conjectures here are true.} that each of us find particularly interesting or intriguing, somewhat in the same style as the first author's --now a decade old-- list of forty-two open problems~\cite{Afonso10L42P}. 

The blog was created and is currently maintained by Afonso S. Bandeira, Daniil Dmitriev, Anastasia Kireeva, Antoine Maillard, Chiara Meroni, Petar Nizic-Nikolac, Kevin Lucca, and Almut R\"{o}dder at ETH Z\"urich. Each blog entry is generally written by one author and the author list of this manuscript is the union of the entry authors in it.

Given the nature of this material, the writing of this manuscript is more informal than a typical academic text\footnote{If you would like to refer to an open question or blog entry, we encourage you to refer to this manuscript instead, since it is a more stable reference.}. Nevertheless, we hope it is useful and inspires thought on these questions. Happy solving!

The numbering of entries and problems follows the blog and is cumulative with~\cite{randomstrasse101-problems2024}.

\tableofcontents

\section{Tensor Concentration Inequalities (KL)}

\medskip \noindent
Consider symmetric deterministic tensors $T_1, \ldots, T_n \in (\mathbb{R}^d)^{\otimes r}$ and let $g_1, \ldots, g_n$ be i.i.d. standard gaussian random variables. We are interested in determining upper bounds for the following quantity:

\begin{equation}\label{eq:expmaxmodel}
    \mathbb{E} \left ( \, \max_{||x||_p \leq 1} \, \left| \sum_{i=1}^n g_i \langle T_i, x^{\otimes r}  \rangle \right| \right)  = \mathbb{E} \left| \left| \sum_{i=1}^n g_i  T_i \right| \right|_{ \mathcal{I}_p}
\end{equation}

\medskip \noindent
Here $ ||x||_p = (\sum_{i=1}^d |x_i|^p)^{\frac1p}$ denotes the $\ell_p$ norm and $\langle T_i, x^{\otimes r} \rangle = \sum_{k_1, \ldots, k_r =1}^d (T_i)_{k_1, \ldots, k_r } x_{k_1} \cdots x_{k_r}$. Moreover, $\left| \left| T_i \right| \right|_{ \mathcal{I}_p} = \max_{||x||_p \leq 1} \left|  \langle T_i, x^{\otimes r}  \rangle \right|$ is the symmetric injective $\ell_p$ norm.

\medskip
\medskip \noindent
Perhaps surprisingly, this question forms an overlap between coding theory, dispersive partial differential equations, tensor principal component analysis, banach space geometry and gaussian process theory (see~\cite{bandeira2024geometric}). Particularly tantalizing to me is the following conjecture:

\begin{conjecture}[Type-$2$ constant of Tensors]\label{conj:type2nck}
    Let $p \geq 2$, then
    \begin{equation}\label{eq:tensorAWNCK}
    \mathbb{E} \left| \left| \sum_{i=1}^n g_i  T_i \right| \right|_{ \mathcal{I}_p} \leq  \, \tilde{\mathcal{O}}_{r,p} \left( d^{\frac12 - \frac1p } \sqrt{\sum_{i=1}^n \left| \left|  T_i \right| \right|_{ \mathcal{I}_p}^2 }\right).
    \end{equation}
    (Here $\tilde{\mathcal{O}}_{r,p}$ hides multiplicative constants depending on $r$ and $p$ as well as polylog factors in $d,n$.)
\end{conjecture}

\medskip \noindent
To get a better sense of what this conjecture is saying, let us focus on the case $r = p = 2$, which is the case where we deal with symmetric matrices $M_i$ and their euclidean operator norm. To the best of my knowledge, the first result of this type had first been established by Tomczak-Jaegermann in~\cite{jaegermann74moduli}, and a more concrete version was shown by Ahlswede and Winter in~\cite{ahlswedewinter02}, which says that for some universal constant $C>0$ we have

\begin{equation}\label{eq:AWmatrixbound}
    \mathbb{E} \left| \left| \sum_{i=1}^n g_i  M_i \right| \right|_{ \mathcal{I}_2} \leq C \sqrt{\log (d+1)} \sqrt{\sum_{i=1}^n \left| \left|  M_i \right| \right|_{ \mathcal{I}_2}^2 }.
\end{equation}

\medskip \noindent
In the scalar case where $d=1$, this statement is of course trivial, as the term on the right side is simply the standard deviation of the gaussian variable on the left multiplied by some constant. The interesting part here is that this inequality tells us that even for much larger $d$, matrices behave similarly to scalars in the sense of this square-root cancellation given by the variance (I consider log factors for now to be a small price to pay). Indeed, the proofs of this inequality mimic techniques used for scalar random variables. The success of this approach lies in the fact that the spectral norm of a symmetric matrix $M$ is well-approximated by traces of higher powers of $M$

\begin{equation}\label{eq:matrixmomentmethod}
    || M || \leq \operatorname{Tr}(M^{2k})^{\frac{1}{2k}} \leq d^{\frac{1}{2k}} || M ||.
\end{equation}

\medskip \noindent
Computing expectations of traces is more manageable and indeed, it is possible to refine \Cref{eq:AWmatrixbound} to arrive at the celebrated noncommutative Khintchine inequality by Lust-Piquard and Pisier~\cite{lustkhinchine86,LPkhintchine91}. Since then numerous generalizations and refinements have appeared~\cite{oli10,Tro12,BBvH-Free,universalityBRvH24}, which have proven to be tremendously useful in applications. 

\medskip \noindent
If we return to the case for general $r,p \geq 2$, we unfortunately lose access to \Cref{eq:matrixmomentmethod} and it is generally NP-hard to approximate injective norms, so one has to come up with a different way to approach the general problem. Observe that \Cref{eq:expmaxmodel} is the supremum of a gaussian process, so its expectation is controlled by a specific distance function on $\mathbb{B}_p^d = \{ x \in \mathbb{R}^d \, \colon \,||x||_p \leq 1\}$ defined as follows:
\begin{equation}
    D(x,y) = \mathbb{E} \left (  \left| \sum_{i=1}^n g_i \langle T_i, x^{\otimes r}-  y^{\otimes r}\rangle \right|^2 \right) 
\end{equation}
If one has decent estimates for the smallest number of balls $\mathcal{N}(\mathbb{B}_p^d,D, \varepsilon)$ with respect to the distance $D$ to cover $\mathbb{B}_p^d$, then one automatically gets a good estimate for the expected injective norm, since Dudley's entropy integral provides the following inequality:

\begin{equation}
    \frac{1}{C \log(n)} \int_{0}^\infty \sqrt{\log \mathcal{N}(\mathbb{B}_p^d,D, \varepsilon)} \, d \varepsilon \, \, \leq  \, \, \mathbb{E} \left| \left| \sum_{i=1}^n g_i  T_i \right| \right|_{ \mathcal{I}_p} \leq C \int_{0}^\infty \sqrt{\log \mathcal{N}(\mathbb{B}_p^d,D, \varepsilon)} \, d \varepsilon
\end{equation}

\medskip \noindent
Unfortunately, controlling $\mathcal{N}(\mathbb{B}_p^d,D, \varepsilon)$ seems to still be a remarkably challenging task. In the case $p=2$ Latała provides a covering number bound for this space which is then combined with an intricate chaining-style argument to get a dimension-free  bound for the injective $\ell_2$ norm of random tensors that was used to get (optimal!) moment estimates for gaussian chaoses~\cite{Lat06}. We (independently) extended this covering number estimate for $p \geq 2$ in~\cite{bandeira2024geometric} to settle \Cref{conj:type2nck} in the regime $p \geq 2r$. Sadly, a volumetric barrier prevents us from proving \Cref{eq:tensorAWNCK} for $p<2r$, which would consist of very interesting cases for further applications (see~\cite{bandeira2024geometric}). Even proving it in the known case $r=p=2$ using purely geometric techniques is, to the best of my knowledge, still open.

\medskip \noindent
To keep the post at a reasonable length, I will refer to~\cite{bandeira2024geometric} for more concrete motivations for this problem and further interesting open questions. 

\medskip \noindent
Lastly, let me point out some interesting related work about norms of nonhomogeneous random tensors (tensors with very irregular randomness patterns) outside of the setting $r=p=2$. Quite recently Aden-Ali removed a logarithmic factor and a constant dependency on $p$ in~\cite{AdenAli25} from our bound in~\cite{bandeira2024geometric} by hammering the problem with the so-called PAC-Bayesian Lemma, which also provides a slick proof of Latała's chaos moment estimate. Boedihardjo proved in~\cite{Boedihardjoinjective24} a remarkably sharp inequality for tensors with independent entries when $p=2$, which is often tight up to a constant. In the independent entry case for matrices people have gone beyond $p=2$ to prove stunningly precise norm estimates~\cite{APSS24normsofstruct,LatalaS25}. The independent entry world seems to allow a fight against logarithmic factors to get dimension-free estimates, while in the general case it seems that we are not even close to proving a crude bound of the form \Cref{eq:tensorAWNCK}. For the case $p=2$ it is also possible to adapt Rudelson's argument in~\cite{Rud96chaining} to get the analogous bound for rank-$1$ tensors (though just as in the independent entry case, it is much more straightforward to prove a bound with suboptimal logs). In~\cite{MASRI2005255} the authors use Hardy space theory to bound the injective $\ell_2$ norm of random hankel tensors and they get in fact a two-sided estimate with the right logarithmic factor. In their paper the tensor has unit-variance rademacher entries, but it can be slightly generalized to allow for other coefficients. In the case $p= \infty$ people have used Kikuchi matrix techniques~\cite{JMBipartiteKikuchi24,BHKLKikuchi24} to get bounds for matching tensors, which tend to have nastier covariance structures, as unlike in the independent entry case, the tensors $T_i$ are not sparse. All of the mentioned examples are consistent with \Cref{conj:type2nck} (and in fact many of them completely outclass it), unfortunately it still seems to be far out of (my) reach.

\section{The Lov\'asz number of random circulant graphs (DD)}

For a graph \(G = (V, E)\), the clique number \(\omega(G)\) and the chromatic number \(\chi(G)\) are fundamental properties studied in combinatorics and graph theory. It is known that computing these quantities is, in general, NP-hard. However, in certain cases, it is possible to efficiently provide upper and lower bounds. As an example, for a \(d\)-regular graph \(G\), it holds that \(\chi(G) \leq d + 1\). In~\cite{lovasz1979shannon}, Lov\'asz proposed a function \(\vartheta(G)\), which was later named the Lov\'asz number or the Lov\'asz theta function, such that it can be efficiently computed and
\begin{equation}
\label{eq:lovasz_omega_chrom}
    \omega(G) \leq \vartheta(\overline{G}) \leq \chi(G),
\end{equation}
where \(\overline{G}\) denotes the complement of the graph \(G\).

The initial motivation of Lov\'asz in~\cite{lovasz1979shannon} was to bound the \emph{Shannon capacity} of a graph. Consider the following setting: label the vertices of a graph \(G\) with distinct letters, and say that two letters can be confused if their respective vertices are adjacent. Similarly, two sequences of length \(k \geq 1\) can be confused, if for each \(1 \leq i \leq k\), the vertices of the \(i\)-th letters are the same or adjacent. Let \(m_k\) be the maximum number of sequences of length \(k\), such that no two sequences can be confused. Clearly, \(m_1 = \alpha(G)\), where \(\alpha(G)\) is the independence number of \(G\). The Shannon capacity of the graph \(G\) is denoted by \(\Theta(G)\) and is defined by 
\begin{equation}
    \Theta(G) := \lim_{k \to \infty} (m_k)^{1/k}.
\end{equation}
Surprisingly, even for simple graphs, such as the cycle on \(7\) vertices \(C_7\), the precise value of \(\Theta(C_7)\) is unknown. The independence number and the Lov\'asz number provide lower and upper bounds on the Shannon capacity, respectively:

\begin{equation}
    \alpha(G) \leq \Theta(G) \leq \vartheta(G).
\end{equation}

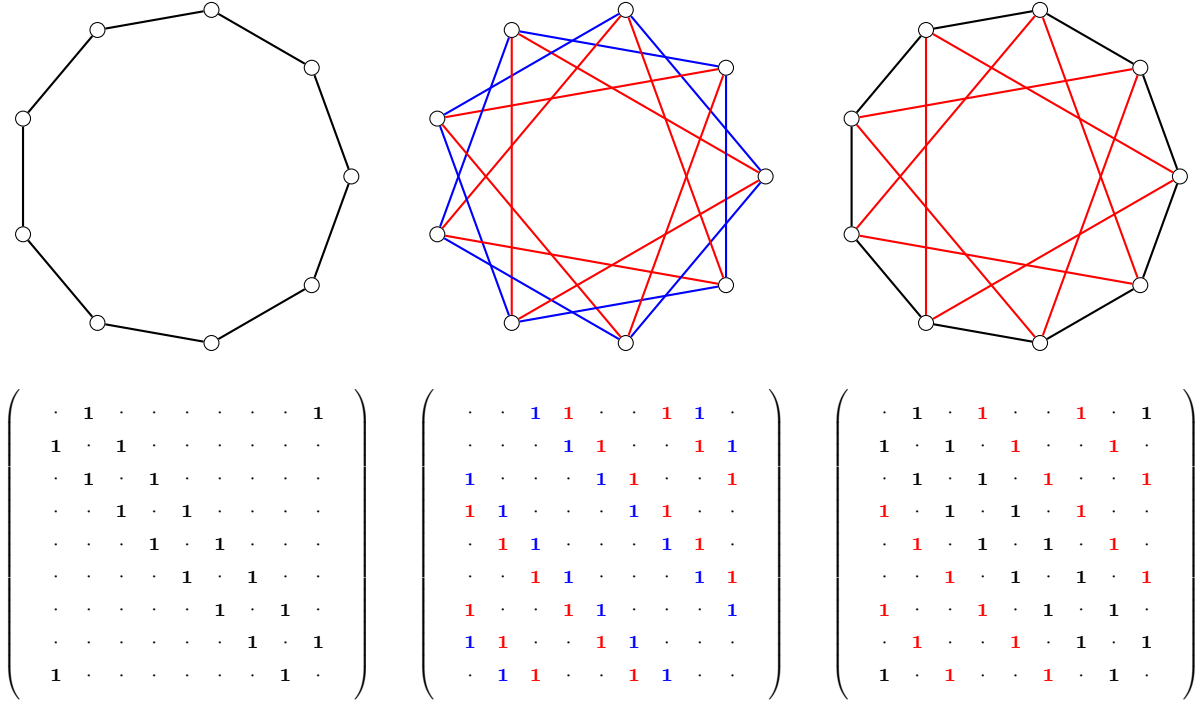
\begin{figure}[h]
\centering
\begin{tabular}{ccc}
\begin{tikzpicture}[scale=1.5]
  \foreach \i in {0,...,8} {
    \node[draw, circle, inner sep=2pt] (v\i) at ({360/9 * \i}:1.5cm) {};
  }
  \foreach \i in {0,...,8} {
    \pgfmathtruncatemacro{\next}{mod(\i+1,9)}
    \draw[thick] (v\i) -- (v\next);
  }
\end{tikzpicture}
&
\begin{tikzpicture}[scale=1.5]
  \foreach \i in {0,...,8} {
    \node[draw, circle, inner sep=2pt] (v\i) at ({360/9 * \i}:1.5cm) {};
  }
  \foreach \i in {0,...,8} {
    \pgfmathtruncatemacro{\nextTwo}{mod(\i+2,9)}
    \draw[blue, thick] (v\i) -- (v\nextTwo);
    \pgfmathtruncatemacro{\nextThree}{mod(\i+3,9)}
    \draw[red, thick] (v\i) -- (v\nextThree);
  }
\end{tikzpicture}
&
\begin{tikzpicture}[scale=1.5]
  \foreach \i in {0,...,8} {
    \node[draw, circle, inner sep=2pt] (v\i) at ({360/9 * \i}:1.5cm) {};
  }
  \foreach \i in {0,...,8} {
    \pgfmathtruncatemacro{\nextOne}{mod(\i+1,9)}
    \draw[thick] (v\i) -- (v\nextOne);
    \pgfmathtruncatemacro{\nextThree}{mod(\i+3,9)}
    \draw[red, thick] (v\i) -- (v\nextThree);
  }
\end{tikzpicture}
\\[1em]
\begin{tikzpicture}[every node/.style={minimum size=0.45cm, anchor=center, font=\scriptsize}]
\matrix (m1) [matrix of nodes,
    left delimiter={(},
    right delimiter={)},
    row sep=-\pgflinewidth, column sep=-\pgflinewidth
]{
  $\cdot$ & $\mathbf{1}$ & $\cdot$ & $\cdot$ & $\cdot$ & $\cdot$ & $\cdot$ & $\cdot$ & $\mathbf{1}$ \\
  $\mathbf{1}$ & $\cdot$ & $\mathbf{1}$ & $\cdot$ & $\cdot$ & $\cdot$ & $\cdot$ & $\cdot$ & $\cdot$ \\
  $\cdot$ & $\mathbf{1}$ & $\cdot$ & $\mathbf{1}$ & $\cdot$ & $\cdot$ & $\cdot$ & $\cdot$ & $\cdot$ \\
  $\cdot$ & $\cdot$ & $\mathbf{1}$ & $\cdot$ & $\mathbf{1}$ & $\cdot$ & $\cdot$ & $\cdot$ & $\cdot$ \\
  $\cdot$ & $\cdot$ & $\cdot$ & $\mathbf{1}$ & $\cdot$ & $\mathbf{1}$ & $\cdot$ & $\cdot$ & $\cdot$ \\
  $\cdot$ & $\cdot$ & $\cdot$ & $\cdot$ & $\mathbf{1}$ & $\cdot$ & $\mathbf{1}$ & $\cdot$ & $\cdot$ \\
  $\cdot$ & $\cdot$ & $\cdot$ & $\cdot$ & $\cdot$ & $\mathbf{1}$ & $\cdot$ & $\mathbf{1}$ & $\cdot$ \\
  $\cdot$ & $\cdot$ & $\cdot$ & $\cdot$ & $\cdot$ & $\cdot$ & $\mathbf{1}$ & $\cdot$ & $\mathbf{1}$ \\
  $\mathbf{1}$ & $\cdot$ & $\cdot$ & $\cdot$ & $\cdot$ & $\cdot$ & $\cdot$ & $\mathbf{1}$ & $\cdot$ \\
};
\end{tikzpicture}
&
\begin{tikzpicture}[every node/.style={minimum size=0.45cm, anchor=center, font=\scriptsize}]
\matrix (m2) [matrix of nodes,
    left delimiter={(},
    right delimiter={)},
    row sep=-\pgflinewidth, column sep=-\pgflinewidth
]{
  $\cdot$ & $\cdot$ & \textcolor{blue}{$\mathbf{1}$} & \textcolor{red}{$\mathbf{1}$} & $\cdot$ & $\cdot$ & \textcolor{red}{$\mathbf{1}$} & \textcolor{blue}{$\mathbf{1}$} & $\cdot$ \\
  $\cdot$ & $\cdot$ & $\cdot$ & \textcolor{blue}{$\mathbf{1}$} & \textcolor{red}{$\mathbf{1}$} & $\cdot$ & $\cdot$ & \textcolor{red}{$\mathbf{1}$} & \textcolor{blue}{$\mathbf{1}$} \\
  \textcolor{blue}{$\mathbf{1}$} & $\cdot$ & $\cdot$ & $\cdot$ & \textcolor{blue}{$\mathbf{1}$} & \textcolor{red}{$\mathbf{1}$} & $\cdot$ & $\cdot$ & \textcolor{red}{$\mathbf{1}$} \\
  \textcolor{red}{$\mathbf{1}$} & \textcolor{blue}{$\mathbf{1}$} & $\cdot$ & $\cdot$ & $\cdot$ & \textcolor{blue}{$\mathbf{1}$} & \textcolor{red}{$\mathbf{1}$} & $\cdot$ & $\cdot$ \\
  $\cdot$ & \textcolor{red}{$\mathbf{1}$} & \textcolor{blue}{$\mathbf{1}$} & $\cdot$ & $\cdot$ & $\cdot$ & \textcolor{blue}{$\mathbf{1}$} & \textcolor{red}{$\mathbf{1}$} & $\cdot$ \\
  $\cdot$ & $\cdot$ & \textcolor{red}{$\mathbf{1}$} & \textcolor{blue}{$\mathbf{1}$} & $\cdot$ & $\cdot$ & $\cdot$ & \textcolor{blue}{$\mathbf{1}$} & \textcolor{red}{$\mathbf{1}$} \\
  \textcolor{red}{$\mathbf{1}$} & $\cdot$ & $\cdot$ & \textcolor{red}{$\mathbf{1}$} & \textcolor{blue}{$\mathbf{1}$} & $\cdot$ & $\cdot$ & $\cdot$ & \textcolor{blue}{$\mathbf{1}$} \\
  \textcolor{blue}{$\mathbf{1}$} & \textcolor{red}{$\mathbf{1}$} & $\cdot$ & $\cdot$ & \textcolor{red}{$\mathbf{1}$} & \textcolor{blue}{$\mathbf{1}$} & $\cdot$ & $\cdot$ & $\cdot$ \\
  $\cdot$ & \textcolor{blue}{$\mathbf{1}$} & \textcolor{red}{$\mathbf{1}$} & $\cdot$ & $\cdot$ & \textcolor{red}{$\mathbf{1}$} & \textcolor{blue}{$\mathbf{1}$} & $\cdot$ & $\cdot$ \\
};
\end{tikzpicture}
&
\begin{tikzpicture}[every node/.style={minimum size=0.45cm, anchor=center, font=\scriptsize}]
\matrix (m3) [matrix of nodes,
    left delimiter={(},
    right delimiter={)},
    row sep=-\pgflinewidth, column sep=-\pgflinewidth
]{
  $\cdot$ & $\mathbf{1}$ & $\cdot$ & \textcolor{red}{$\mathbf{1}$} & $\cdot$ & $\cdot$ & \textcolor{red}{$\mathbf{1}$} & $\cdot$ & $\mathbf{1}$ \\
  $\mathbf{1}$ & $\cdot$ & $\mathbf{1}$ & $\cdot$ & \textcolor{red}{$\mathbf{1}$} & $\cdot$ & $\cdot$ & \textcolor{red}{$\mathbf{1}$} & $\cdot$ \\
  $\cdot$ & $\mathbf{1}$ & $\cdot$ & $\mathbf{1}$ & $\cdot$ & \textcolor{red}{$\mathbf{1}$} & $\cdot$ & $\cdot$ & \textcolor{red}{$\mathbf{1}$} \\
  \textcolor{red}{$\mathbf{1}$} & $\cdot$ & $\mathbf{1}$ & $\cdot$ & $\mathbf{1}$ & $\cdot$ & \textcolor{red}{$\mathbf{1}$} & $\cdot$ & $\cdot$ \\
  $\cdot$ & \textcolor{red}{$\mathbf{1}$} & $\cdot$ & $\mathbf{1}$ & $\cdot$ & $\mathbf{1}$ & $\cdot$ & \textcolor{red}{$\mathbf{1}$} & $\cdot$ \\
  $\cdot$ & $\cdot$ & \textcolor{red}{$\mathbf{1}$} & $\cdot$ & $\mathbf{1}$ & $\cdot$ & $\mathbf{1}$ & $\cdot$ & \textcolor{red}{$\mathbf{1}$} \\
  \textcolor{red}{$\mathbf{1}$} & $\cdot$ & $\cdot$ & \textcolor{red}{$\mathbf{1}$} & $\cdot$ & $\mathbf{1}$ & $\cdot$ & $\mathbf{1}$ & $\cdot$ \\
  $\cdot$ & \textcolor{red}{$\mathbf{1}$} & $\cdot$ & $\cdot$ & \textcolor{red}{$\mathbf{1}$} & $\cdot$ & $\mathbf{1}$ & $\cdot$ & $\mathbf{1}$ \\
  $\mathbf{1}$ & $\cdot$ & \textcolor{red}{$\mathbf{1}$} & $\cdot$ & $\cdot$ & \textcolor{red}{$\mathbf{1}$} & $\cdot$ & $\mathbf{1}$ & $\cdot$ \\
};;
\end{tikzpicture}
\\[0.5em]
\end{tabular}
\caption{Top row: Circulant graphs on 9 vertices; Bottom row: Their corresponding adjacency matrices (with 0's represented by dots).}
\label{fig:circ_graphs}
\end{figure}

In the following, we discuss how to bound the L\'ovasz number for a certain family of graphs.
Recall that the independence number, which is equal to the largest size of an independent set, can be defined as follows (\(E(G)\) is the set of edges):

\begin{equation}
\label{def:alpha}
    \alpha(G) := \max_{v \in \{0, 1\}^n}
    \Big\{\sum_{i} v_i \text{, such that } v_i v_j = 0\text{ for all }(i, j) \in E(G), \Big\}.
\end{equation}

Lov\'asz number is a natural semidefinite programming (SDP) relaxation for~\Cref{def:alpha}:

\begin{equation}
\label{def:ltn}
    \vartheta(G) := \max_{X \in \mathbb{R}^{n \times n}}
    \Big\{\sum_{i, j} X_{ij} \text{, such that } X \succeq 0, \mathrm{Tr}\,X = 1, \text{ and } X_{ij} = 0\text{ for all }(i, j) \in E(G)\Big\}.
\end{equation}

An interesting question is to characterize the Lov\'asz number for random graphs. 
A well-known family of random graphs is the Erd\H{o}s-R\'enyi random graphs \(G(n, p)\). Here, \(n\) denotes the number of vertices, and \(0 \leq p \leq 1\) is the probability parameter. To sample a graph from \(G(n, p)\), for each pair of vertices \(i \neq j \in [n]\), one adds an edge independently with probability \(p\). The random graphs \(G(n, p)\) are ubiquitous in mathematics; for instance, they are used in combinatorics to lower bound \emph{Ramsey} numbers. There are several results on the Lov\'asz number for Erd\H{o}s-R\'enyi random graphs \(G(n, p)\) (see~\cite{juhasz1982asymptotic,coja2005lovasz,aroranote}). However, the precise asymptotics of the Lov\'asz number for \(G(n, 1/2)\) is still open.

\begin{conjecture}
    Let \(G\) be distributed as \(G(n, 1/2)\). Then,
    \begin{equation}
        \mathbb{E} \vartheta(G) = (1 + o(1))\sqrt{n}.
    \end{equation}
\end{conjecture}

Lov\'asz number exhibit a beautiful property~\cite{lovasz1979shannon}: for any graph \(G\)
\begin{equation}
\label{eq:lovasz_complement}
    \vartheta(G) \vartheta(\overline{G}) \geq n.
\end{equation}
Note the connection to the inequality \(w(G)\chi(\overline{G}) \geq n\) via~\Cref{eq:lovasz_omega_chrom}. 
For vertex-transitive graphs \(G\), which include circulant graphs,~\Cref{eq:lovasz_complement} holds with equality.

Circulant graphs~(see~\Cref{fig:circ_graphs}) are examples of highly structured graphs and are defined as Cayley graphs over the cyclic group \(\mathbb{Z}_n\). In other words, their adjacency matrix is circulant: if \((i, j) \in E(G)\), then \((i + 1\ \mathrm{mod}\ n, j+1\ \mathrm{mod}\ n) \in E(G)\). A celebrated example of a circulant graph from number theory is \emph{the Paley graph}. For a prime \(p \equiv 1\ \mathrm{mod}\ 4\), in order to define the Paley graph of order \(p\), consider the vertex set \([p]\) and connect two vertices \(i\) and \(j\) if and only if \(i - j\) is a quadratic residue modulo \(p\). Paley graphs are of great importance in combinatorics and theoretical computer science, since they are believed to behave \emph{pseudorandomly}: they are deterministic graphs which share certain properties with a typical realization of \(G(p, 1/2)\). Current rigorous evidence of this phenomenon relies strongly on results from number theory and additive combinatorics.

Random (dense) circulant graphs can be viewed as a middle step between a 'fully random' \(G(p, 1/2)\) and a 'fully deterministic' Paley graph. To sample a random circulant graph, we only need to sample the neighbors of a single vertex, say vertex \(0\). Since \((0, k) \in E(G)\) implies that \((0, n - k) \in E(G)\), we sample a vector \(x \in \{0, 1\}^{\left\lceil\frac{n - 1}{2}\right\rceil}\) uniformly at random, and \(x_k = 1\) stands for \((0, k) \in E(G)\).

\begin{table}[t]
\centering
\caption{Four equivalent LPs for \(\vartheta(G)\).}
\label{table:4_lps}
\renewcommand{\arraystretch}{1.05}
\begin{tabular}{@{}p{0.47\linewidth} p{0.47\linewidth}@{}}
\multicolumn{2}{c}{\textbf{'time' domain}} \\[3pt]
\toprule
\textbf{Primal} & \textbf{Dual}\\
\midrule
\(
\begin{aligned}
\max_{x \in \mathbf{R}^n} \ &\sum_{i} x_i\\
\text{s.t. } &x_k = x_{n - k} \\
&\quad \text{for all } k \in [n] \setminus \{0\}, \\
&x_0 = 1,\ Fx \ge \mathbf{0} ,\\
&x_k=0 \\
&\quad \text{for all }(0, k) \in E(G).
\end{aligned}
\)
&
\(
\begin{aligned}
\min_{z \in \mathbf{R}^n} &1 + \sum_i z_i \\
\text{s.t. } &z_k = z_{n - k} \\
&\quad \text{for all } k \in [n] \setminus \{0\}, \\
&z \geq \mathbf{0}, \\
&\langle z, f_k \rangle = -1 \\
&\quad \text{for all } (0, k) \in E(\overline{G}).
\end{aligned}
\)
\\
\bottomrule
\\[-0.2em]
\multicolumn{2}{c}{\textbf{'frequency' domain}} \\[3pt]
\toprule
\textbf{Primal} & \textbf{Dual}\\
\midrule
\(
\begin{aligned}
\max_{y \in \mathbf{R}^n}\ &n y_0\\
\text{s.t. } &y_k = y_{n - k} \\
&\quad \text{for all } k \in [n] \setminus \{0\},  \\
&\lVert y \rVert _1 = 1,\ y \geq \mathbf{0},\\
&\langle y, f_k \rangle = 0 \\
&\quad \text{for all } (0, k) \in E(G).
\end{aligned}
\)
&
\(
\begin{aligned}
\min_{t \in \mathbf{R}^n} &1 + n t_0\\
\text{s.t. } &t_k = t_{n - k} \\
&\quad \text{for all } k \in [n] \setminus \{0\},\\
&Ft \geq \mathbf{0}, \\
&t_k = -1 / n \\
&\quad \text{for all } (0, k) \in E(\overline{G}).
\end{aligned}
\)
\\
\bottomrule
\end{tabular}
\end{table}

Let \(F \in \mathbb{C}^{n \times n}\) be the discrete Fourier transform matrix: \(F_{jk} = \exp(-2 \pi i jk/n)\) for \(j, k \in [n]\). Here and in the following, we index vectors and matrices by \([n] := \{0, 1, \ldots, n - 1\}\). Let \(f_k\) denote the \(k\)-th row of \(F\).
\cite{magsino2019linear}~notes that, since circulant matrices are diagonalizable by \(F\), the Lov\'asz number for circulant graphs can be written as a linear program. In~\Cref{table:4_lps}, we present four equivalent linear programs arising from strong duality (see, e.g.,~\cite{bazaraa2011linear}) and switching between the 'time' and the 'frequency' domains. For the latter, we perform the change of variables, \(y := Fx\) and \(t := Fz\), respectively. 

We conjecture, based on numerical observations, that the Lov\'asz number for random dense circulant graphs behaves similarly to that for \(G(n, 1/2)\). This motivates the following question.

\begin{conjecture}
\label{prob:lovasz_circ}
    Let \(G\) be distributed as a random dense circulant graph. Then,
    \begin{equation}
        \mathbb{E} \vartheta(G) = (1 + o(1))\sqrt{n}.
    \end{equation}
\end{conjecture}

In~\cite{bandeira2025lov}, we provide a partial result towards resolving~\Cref{prob:lovasz_circ}, with a precise lower bound and an upper bound \(O(\sqrt{n \log \log n})\). Our techniques rely on studying a modified version of a linear program from~\Cref{table:4_lps} in primal form and frequency domain, and using the restricted isometry property (RIP) of the subsampled discrete Fourier transform matrix.

\section{Injectivity and Stability of Phase Retrieval (ASB)}

Throughout this post $\mathbb{K}$ will stand for either $\mathbb{R}$ or $\mathbb{C}$.

Given $A\in\mathbb{K}^{N\times M}$, the Phase Retrieval Problem aims to recover a vector $x\in\mathbb{K}^M$ from $|Ax|$, where $|\cdot|$ is the entrywise modulus. In other words, one has access to the modulus of $N$ linear measurements, and the goal is to recover $x\in\mathbb{K}^M$. Let $\mathbb{T}_\mathbb{K}$ denote the unit torus in $\mathbb{K}$ (corresponding to the unit circle in $\mathbb{C}$ or $\{\pm1\}$ in $\mathbb{R}$). One can only hope to recover $x$ modulo $\mathbb{T}_\mathbb{K}$ (meaning up to a global phase shift in $\mathbb{C}$ or a global sign flip in $\mathbb{R}$).

We say $A$ is injective for phase retrieval (for $\mathbb{K}$) if the map 
\[
\mathcal{A}:x \bmod \mathbb{T}_\mathbb{K} \mapsto |Ax|,
\]
is injective.

In 2013 a few collaborators and I~\cite{Bandeira_etal_Savingphase} conjectured that $N\geq 4M-4$ linear measurements were necessary for injectivity over $\CC$, and that generic $4M-4$ measurements were injective. While the second part of this conjecture was proven in~\cite{conca2015generico4Mm4}, the first was shown to be false when Cynthia Vinzant found an injective set of $11$ linear measurements in $\mathbb{C}^4$~\cite{vinzant2015small}.\footnote{A prize was awarded for this finding:~\url{https://dustingmixon.wordpress.com/2015/07/08/conjectures-from-sampta/}.} Cynthia stated a refined version of the conjecture that still remains open (see~\url{https://dustingmixon.wordpress.com/2015/07/08/conjectures-from-sampta/})

\begin{conjecture}
Let $N=4M-5$ and draw $A\in\CC^{N\times M}$ at random (here the important thing is for $ \mathrm{im}(A) $ to be drawn uniformly from the Grassmannian of $ M $-dimensional subspaces of $ \mathbb{C}^{4M-5}$, it can e.g. be with iid gaussian entries). Let $p_M$ be the probability that the mapping $ x \bmod \mathbb{T} \mapsto |Ax|^2 $ is injective.

\begin{itemize}
    \item[(a)] $ p_M < 1 $ for all $ M $.
    \item[(b)] $ \displaystyle\lim_{M \to \infty} p_M = 0 $.
\end{itemize}

\end{conjecture}

Over $\mathbb{R}$ the question of injectivity is significantly easier and it is known that $N\geq 2M-1$ is necessary that that generic $2M-1$ measurements are injective. A crucial step in the proof of this statement is the fact that any partition of the rows of $A^{(2M-1)\times M}$ into two sets yields one set that spans $\mathbb{R}^M$. This is known as the complement property (see, e.g.,~\cite{Bandeira_etal_Savingphase}). 

A natural question is whether such phaseless measurements allow for stable reconstruction of $x$. In this direction, Balan and Wang~\cite{balan2013invertibility} showed that stability of the map $\mathcal{A}$ is related to the condition number of subsets of rows of $A$ that span $\mathbb{R}^M$.

\begin{definition}
Given $A\in\mathbb{R}^{N\times M}$ let $\omega(A)$ be defined as
\begin{equation}
\omega(A) = \min_{S\subseteq [N]:\, \mathrm{rank}(A_{S^c})<n} \sigma_M(A_S),
\end{equation}
where $A_S$ denotes the submatrix made of the rows indexed by $S$ and $\sigma_M$ is the $M$-th singular value (which is the smallest in non-degenerate situations).
\end{definition}
The complement property says that $\mathcal{A}$ is injective if and only if $\omega(A)>0$.~\cite{balan2013invertibility} shows that stability is tightly connected to the inverse of this quantity. They also make the following intriguing conjecture.

\begin{conjecture}
    There exist universal constants $C>0$ and $0<\beta<1$ such that, for any $M>1$ and $A\in\mathbb{R}^{2M-1\times M}$ a matrix for which any subset of $M$ rows spans $\mathbb{R}^M$,
    \[
    \omega(A) \leq C \max_{k\in[N]}\|A_k\|\beta^M,
    \]
    where $A_k$ is the $k$-th row of $A$.
\end{conjecture}

This also motivates the following question (for which I do not know the state of the art, but a quick search did not reveal that the answer is already known).

\begin{openproblem}
    What is the value of $\omega(A)$ for $A$ with iid standard gaussian entries? Does it satisfy the bound in the conjecture? If so, what is $\beta$?
\end{openproblem}

\textbf{Acknowledgments:} I would like to thank Mitchell Taylor from bringing the Balan-Wang conjecture to my attention.

\section{Mutually Unbiased Bases, ETFs, and Zauner's Conjecture (ASB)}

Finite frame theory is the source of many beautiful questions. A remarkable example concerns the existence of Mutually Unbiased Bases,  objects of interest in Quantum Physics (see~\cite{McNulty2024-MUBs} for a review).

\begin{definition}[Mutually Unbiased Bases]
    Two orthonormal bases $\{v_1, \dots, v_d\}, \{w_1, \dots, w_d\} \subset \mathbb{C}^d$ are called \emph{mutually unbiased} if $|\langle v_i, w_j \rangle| = 1 / \sqrt{d}$ for all $i, j$. Let $\mathrm{MUB}(d)$ denote the largest possible number of (simultaneously) mutually unbiased bases (MUBs) in $\mathbb{C}^d$.
\end{definition}

In~\cite{Bandyopadhyay2002-MUBs} it is shown that $\mathrm{MUB}(d)\leq d+1$ and that equality is achieved when $d$ is a prime power. For $d$ not a prime power, the problem of determining $\mathrm{MUB}(d)$ is still open. The smallest open instance is particularly well known for being a tantalizing open problem that has remained open (see Open Problem 6.2 in~\cite{Afonso10L42P}).

\begin{conjecture}
    Let $\mathrm{MUB}(d)$ denote the largest possible number of (simultaneously) mutually unbiased bases (MUBs) in $\mathbb{C}^d$. We have $\mathrm{MUB}(6)<7$.
\end{conjecture} 

It is known that $\mathrm{MUB}(6)\geq 3$ (see, e.g.,~\cite{BDK-2022-DualBoundsMUB}). Note that for $\RR^d$ the existence of two mutually unbiased bases is equivalent to the existence of Hadamard Matrices, Conjeture 14 in this blog~\cite{randomstrasse101-problems2024}.

An interesting approach to this problem is to try to build Sum-of-Squares proofs of $\mathrm{MUB}(6)<7$. In~\cite{BDK-2022-DualBoundsMUB} it was shown that a certain relaxation related to degree-2 Sum-of-Squares cannot prove that $M(6) < 7$. It is however unclear whether higher degree levels of Sum-of-Squares can provide such a proof.

\begin{openproblem}
    Is there a Sum-of-Squares degree 4 proof that there are no $7$ Mutually Unbiased Bases in $\mathbb{C}^6$?
\end{openproblem}

Note that this is a statement about the existence, or not, of $7\times 6$ vectors $v_{i}^{(k)}\in \mathbb{C}^6$ for $i=1,\dots,6$ and $k=1,\dots,7$ satisfying
\[
    \left| \left\langle v_{i}^{(k)},v_{j}^{(\ell)}\right\rangle \right|^2 = \left\{ \begin{array}{ccl}
    1 & \text{if} & k= \ell \text{ and } i= j\\
    0 & \text{if} & k= \ell \text{ and } i\neq j\\
    \frac{1}{\sqrt{6}} & \text{if} & k\neq \ell.
    \end{array} \right.
    \]
These can for example be encoded with quartic relations on the real and imaginary parts of the vectors components (since $\left(v^\ast u\right)\left(u^\ast v\right) = |\langle u, v \rangle|^2$). An alternative encoding can be written by noting that the constraints for $k=\ell$ are quadratic and that, since $1/\sqrt{6}$ is the smallest possible value for the largest inner-product between vectors in different basis, one can rewrite the $k\neq \ell$ constraints as $| \langle v_{i}^{(k)},v_{j}^{(\ell)}\rangle |\leq 1/\sqrt{6}$ which are semidefinite constraints (on $2\times 2$ Hermitian matrices) involving a quadratic quantity: $\left[\begin{array}{cc}1/\sqrt{6} & \langle v_{i}^{(k)},v_{j}^{(\ell)}\rangle \\ \langle v_{j}^{(\ell)},v_{i}^{(k)}\rangle & 1/\sqrt{6}\end{array}\right]\succeq 0$.

Another very interesting object in this field is an Equiangular Tight Frame (ETF). In $\mathbb{C}^d$ a frame is a set of spanning vectors $\phi_1,\dots,\phi_n\in \mathbb{C}^d$. A frame is called unit-normed if all vectors have unit-norm and it is called a tight frame if $\sum_{i=1}^n \phi_i\phi_i^\ast$ is a multiple of the $d\times d$ identity matrix. In general the smallest eigenvalue of this matrix is called the lower frame bound, and the largest the upper frame bound, a tight frame corresponds to a frame for which these match (see \S 11 in~\cite{Bandeira2025MathofSNL}). The worst-case coherence $\mu$ of a unit-norm frame $\phi_1,\dots,\phi_n$ is given by $\mu:=\max_{i\neq j} |\langle \phi_i,\phi_j\rangle|$ (note that Mutually Unbiased Bases, described above, have worst-case coherence $1/\sqrt{d}$). Equiangular tight frames are the unit-norm frames with the smallest worst-case coherence among frames with the same size and dimension (see \S 12 in~\cite{Bandeira2025MathofSNL}).

\begin{definition}[Equiangular Tight Frame (ETF) in $\mathbb{C}^d$]
    A frame $\phi_1,\dots,\phi_n\in \mathbb{C}^d$ is an ETF if:
    \begin{enumerate}
        \item It is unit-normed: $\|\phi_i\|=1$ for all $i=1,\dots,n$.
        \item It is a tight frame: $\sum_{i=1}^n \phi_i\phi_i^\ast$ is a multiple of the $d\times d$ identity matrix.
        \item It is equiangular: there exists $\mu\geq 0$ such that $|\langle \phi_i,\phi_j\rangle|=\mu$ for all $i\neq j$. 
    \end{enumerate}
    Note that $\mu$ needs to match the Welch bound, meaning that $\mu = \sqrt{(N-d)/(N-1)}/\sqrt{d}$.
\end{definition}

With a tensoring argument (see \S 12 in~\cite{Bandeira2025MathofSNL}) one can show that $n\leq d^2$. ETFs in $\mathbb{C}^d$ with $n=d^2$ vectors are particularly important in Quantum Physics and are known as \texttt{SIC-POVM}s (Symmetric, Informationally Complete, Positive Operator-Valued Measures). It is a well known conjecture that they exist for all dimensions $d$, known as Zauner's Conjecture~\cite{Zauner1999}.

\begin{conjecture}[Zauner's Conjecture]
    For each dimension $d$, there exists an ETF in $\mathbb{C}^d$ with $n=d^2$ vectors.
\end{conjecture} 

There has been fascinating recent progress~\cite{ApplebyFlammiaKopp2025StarkSIC}: a proposed construction (corresponding to an orbit of the Weyl–Heisenberg group acting on $\mathbb{C}^d$) has been conditionally proven to indeed be a \texttt{SIC-POVM} assuming a couple of conjectures in Number Theory, examples of so-called Stark Conjectures. This shows that these conjectures imply Zauner's Conjecture. A unconditional proof of Zauner's Conjecture remains unavailable.

\section{On the clique number of the Paley Graph (ASB)}

Given a prime number $p$ such that $p \equiv 1 \pmod{4}$ the Paley graph on $p$ nodes $G_p$ is the graph where nodes $i$ and $j$ are connected if $i-j$ is a quadratic residue modulo $p$. This graph is a regular graph with degree $\frac{p-1}{2}$ and it is believed to have many pseudo-random properties (mimicking in many ways a uniformly random regular graph of the same degree, or even a $G(n,1/2)$). See~\cite{Kunisky-2023-SpectralPseudorandomnessPaleyClique} for a nice overview aligned with the topic of this post.

Let $\omega(G_p)$ denote the clique number of $G_p$ i.e. the size of the largest set of nodes that are all pairwise connected by an edge. 
Since $G_p$ is isomorphic to its complement this is the same as $\alpha(G_p)$ the independence number of $G_p$ (the size of the largest set of nodes that have no edges between them). We will focus on the clique number of ease of exposition. A random graph of the same degree has logarithmic clique number, which motivates the following well known conjecture (see, e.g., Open Problem 8.4 in~\cite{Afonso10L42P})

\begin{conjecture}
    Let $\omega(G_p)$ denote the clique number of the Paley Graph for $p \equiv 1 \pmod{4}$ prime.
    \[
    \omega(G_p)=O(\mathrm{polylog(p)}).
    \]
\end{conjecture} 

An upper bound of $\omega(G_p)\leq \sqrt{p}$ can be shown either by a spectral method or by computing the Lovasz's theta function $\vartheta(\overline{G_p})$ of the complement of $G_p$ (the Lovasz's theta function corresponds to an SDP and is defined on Post 9 of this blog) and recalling that $\omega(G_p)\leq \vartheta(\overline{G_p})$. Recently, Hanson and Petridis~\cite{HP-2021-SumsetsPaleyClique} showed that $\omega(G_p) \leq (1 + o(1))\sqrt{p/2}$. The post is mostly devoted to the question of whether convex relaxations such as the Lovasz's theta function can provide interesting upper bounds on this quantity.

We start with the notion of localization~\cite{MMP-2019-LPCliquesPaley,Kunisky-2023-SpectralPseudorandomnessPaleyClique}: Let $G_{p,1}$ denote the induced subgraph on the vertices adjacent to the node 0, which we call the 1-localization. Because the 1-localization does not depend on the node picked (see~\cite{Kunisky-2023-SpectralPseudorandomnessPaleyClique}) we have $\omega(G_p) \leq \omega(G_{p,1}) + 1 \leq \vartheta(\overline{G_p})+1$. It was observed empirically in \cite{MMP-2019-LPCliquesPaley} that $\vartheta(\overline{G_{p,1}}) \sim \sqrt{p/2}$, which would recover the Hanson-Petridis bound (also, note that $G_{p,1}$ is a circulant graph and so the SDP corresponding to $\vartheta(\overline{G_{p,1}})$ reduces to an LP).

\begin{conjecture}
    $\vartheta(\overline{G_{p,1}}) \sim \sqrt{p/2}$ (for $p \equiv 1 \pmod{4}$ prime) where $G_{p,1}$ is the 1-localization of the Paley Graph described above.
\end{conjecture}

This process can be iterated: Let $G_{p,2}$ denote the induced subgraph on the vertices adjacent to 0 and 1. While the resulting graph is no longer circulant, it is still the case that there is only one 2-localization (up to isomorphism), and so $\omega(G_p) \leq \omega(G_{p,2}) + 2\leq \vartheta(\overline{G_{p,2}})+2$. \cite{Kunisky-2023-SpectralPseudorandomnessPaleyClique} observed empirically that $\omega(G_{p,2}) \sim (\sqrt{1/2} - \epsilon)\sqrt{p}$ (for some $\epsilon>0$), suggesting this approach could improve on the Hanson-Petridis bound. Perhaps $\omega(G_{p,2}) \sim \frac{\sqrt{p}}{2}$). We note that for 3-localizations several graphs would have to be checked (picking different triplets of nodes may yield non-isomorphic subgraphs).

\begin{conjecture}
        For $p \equiv 1 \pmod{4}$ prime and large enough we have $\vartheta(\overline{G_{p,2}}) \leq \frac23 \sqrt{p}$ where $G_{p,2}$ is the 2-localization of the Paley Graph described above.
\end{conjecture}

We note that Conjectures 17 and 18 in this Blog are consistent with these conjectures and the fact that one expects $G_p$ and their localizations to have pseudo-random properties.

There is another approach to to use convex relaxation to potentially improve over upper bounds on the $
\omega(G_p)$. Recall that the Lovasz $\vartheta$ function can be viewed a as degree-2 Sum-of-Squares relaxation of the clique (or independence) number. A fascinating question is whether the degree-4 Sum-of-Squares relaxation significantly improves on the $\sqrt{p}$ bound, perhaps even improving over the state of the art (see~\cite{KY-2022-PaleyDegree4SOS} for a precise construction of this relaxation). There have been several encouraging numerical experiments (for this and other SDPs)~\cite{KY-2022-PaleyDegree4SOS,kobzar2023revisiting,wang2024lower}. We particularly note that the numerical experiments in~\cite{KY-2022-PaleyDegree4SOS} suggest even a polynomial improvement may be possible, although the authors prove that it cannot improve beyond $\Omega(p^{\frac13})$.

\begin{conjecture}
    There exists $\varepsilon>0$ such that the Sum-of-Squares relaxation of degree 4 for the clique number of the Paley graph gives a bound of $O(p^{\frac12-\varepsilon})$.
\end{conjecture}

The Paley Clique Conjecture has been formally related to a conjecture related to sparse recovery known as the Paley ETF conjecture~\cite{Bandeira_ConditionalPaley,Bandeira_RoadRIP}. For $p \equiv 1 \pmod{4}$ prime the Paley ETF (see Definition 21.8 in~\cite{Bandeira2025MathofSNL}) is an Equiangular Tight Frame (as defined in post 11 of this blog) with $n=p+1$ vectors in $\mathbb{C}^{\frac{p+1}{2}}$. We defer the detailed construction to the references above but in a nutshell the matrix whose columns are this frame is built by taking a submatrix of the Discrete Fourier Transform matrix corresponding to rows whose indices are quadratic residues (after an appropriate scaling a canonical basis column is added, see Definition 21.8 in~\cite{Bandeira2025MathofSNL} or~\cite{Bandeira_ConditionalPaley,Bandeira_RoadRIP}). The resulting $\frac{p+1}{2}\times (p+1)$ matrix is conjectured to satisfy the Restricted Isometry Property beyond the so-called square-root bottleneck~\cite{Tao_blog_deterministicRIP,Bandeira_RoadRIP} (see also Open Problem 6.4 in \cite{Afonso10L42P}).

\begin{conjecture}
    The Paley Equiangular Tight Frame construction satisfies the Restricted Isometry Property beyond the square-root bottleneck. In other words, there exists $\varepsilon>0$ such that, for $m \sim p^{1/2+\varepsilon}$, any matrix formed by picking $m$ (distinct) vectors in the Paley ETF has condition number bounded by a constant independent of $p$. 
\end{conjecture} 

\section{The KLS Conjecture and Implications (AR)}

The isoperimetric problem dates back to the ancient Greeks who conjectured that in $\mathbb{R}^2$, the shape with fixed boundary length that  maximizes area is given by a circle.
While this has been proven a long time ago, the same question for higher-dimensional spaces and general measures remains an active area of research.
For a probability measure $\mu$ on $\mathbb{R}^n$, and Euclidean distance $d$, the isoperimetric problem consists in determining the set $S$ that minimizes the boundary measure $\mu^+(S)$ given a fixed volume $\mu(S)$.
\begin{definition}
    The boundary measure of a set $S \subset \mathbb{R}^n$ with respect to the probability measure $\mu$ and Euclidean distance $d$ is given by
    \begin{align*}
        \mu^+(S):= \lim_{\varepsilon \to 0}\frac{\mu(S+\varepsilon B_n)-\mu(S)}{\varepsilon}. 
    \end{align*}
    where $B_n$ denotes the Euclidean Ball in $n$ dimensions.
    The isoperimetric  constant $\psi_\mu$ is the minimal boundary measure among sets $S$ with $\mu(S)\geq \frac{1}{2}$:
    \begin{align*}
        \psi_\mu := \inf_{S \subseteq \mathbb{R}^n} \frac{\mu^+(S)}{\min\{\mu(S), \mu(\mathbb{R}^n\backslash S)\}}. 
    \end{align*}
\end{definition}
Sometimes $\psi_\mu$ is referred to as the Cheeger constant.
For the $n$-dimensional standard Gaussian measure $\gamma_n$, the isoperimetric constant is given by the minimum over half-spaces (first proven in~\cite{sudakov1974extremal}, but several proofs exist nowadays). Notably, this implies lower bounds on the isoperimetric constant independent of the dimension.

On the other hand, a general measure $\mu$ might admit a bottleneck, i.e. a split  dividing $\mathbb{R}^n$ into two or more regions that are only connected by small probability sets (the reader may want to image a gym dumbbell), and in those cases the isoperimetric constant can be arbitrarily small. 
This motivates the following question: Does there exist a class of measures on $\mathbb{R}^n$ where we do not expect bottlenecks, or more precisely, such that the isoperimetric constant is bounded from below by a dimension-independent and universal constant?

Consider the class of log-concave probability distributions. These are the probability measures with densities $f$ such that $\log(f)$ is concave. This includes in particular all uniform measures over convex sets, a situation where we do not expect bottlenecks to appear.
Indeed, Kannan, Lovász and Simonovits~\cite{kannan1995isoperimetric} conjectured in 1995 that for all log-concave measures, the isoperimetric constant is, up to a universal constant, attained by minimizing over half-spaces, like in the Gaussian case:
\begin{conjecture}[KLS Conjecture (\cite{kannan1995isoperimetric})]
There exists a universal constant $C>0$ such that for every log-concave measure $\mu$, the isoperimetric constant satisfies
     \begin{align*}
         \psi_\mu \geq C \inf_{H \subseteq \mathbb{R}^n, H \text{ is halfspace}} \frac{\mu^+(S)}{\min\{\mu(S), \mu(\mathbb{R}^n\backslash S)\}}.
     \end{align*}
     Equivalently, there exists a universal constant $C>0$ such that
    \begin{align*}
         \psi_\mu \geq \frac{C}{\sqrt{\|A\|_\text{op}}},
     \end{align*}
     where $A$ denotes the covariance matrix of the measure $\mu$, $A=\mathbb{E} \left(X- \mathbb{E}X\right)\left(X- \mathbb{E}X\right)^\top$, where $\EE$ denotes expectation of $X\sim\mu$. 
\end{conjecture}
The equivalence is due to log-concavity and shown in Lemma 1 in~\cite{lee2018kannan}. Moreover, by this equivalence we can focus on isotropic distributions, i.e. $\mathbb{E}[X]=0$ and $A= I$. We will denote by $\psi_n$ the infimum over all log-concave isotropic distributions over $\mathbb{R}^n$.

Besides important connections to other conjectures which we will discuss later, a lower bound on the isoperimteric constant has very useful implications:
\begin{itemize}
    \item It is equivalent to a Poincaré inequality which states that there exists a constant $C_P>0$ such that for all smooth functions $f: \mathbb{R}^n \to \mathbb{R}$,
    \begin{align*}
        \mathrm{Var}_{X \sim \mu}(f(X)) \leq C_P \int_{\mathbb{R}^n} \|\nabla f\|_2^2 \,d \mu.
    \end{align*}
    
    \item Moreover, a lower bound on the isoperimetric constant is equivalent to Lipschitz concentration (\cite{gromov1983topological},~\cite{milman2009role}): For $X \sim \mu$ and $f$ being an $L$-Lipschitz function, it holds 
    \begin{align*}
        \mathbb{P}(|f(X)-\mathbb{E}[f(X)]|\geq t) \leq \exp(-\frac{t}{L}\Omega(\psi_\mu)).
    \end{align*}
    \item The ball walk on a convex body is a Markov Chain Monte Carlo algorithm to sample from the uniform distribution over a convex body. If the starting distribution is sufficiently close to the uniform distribution, the KLS conjecture implies that for every convex body the ball walk mixes in $\tilde{O}(n^2)$ steps (up to log factors (\cite{kannan1997random})). This rate is, in general, unimprovable.
\end{itemize}
For a more detailed description and further applications, see the survey by Lee and Vempala~\cite{lee2018kannan}.

In recent years, there has been significant progress on improving the lower bound for the isoperimetric constant uniformly over log-concave distributions. The starting point in~\cite{kannan1995isoperimetric} was the following:
For one-dimensional log-concave distributions, it is well known (and simple) to deduce by concavity that the isoperimetric constant is uniformly bounded. This fact was used to reduce high-dimensional log-concave distributions to one-dimensional distributions through \textit{localization}. This technique was used to prove $\psi_\mu \geq \frac{C}{n^\frac{1}{2}}$ for isotropic distributions in~\cite{kannan1995isoperimetric}.
Besides many important improvements, we briefly summarize the most influential improvements throughout the past fifteen years. Based on the fact that
Gaussians and measures which are \textit{at least as} log-concave as a Gaussians, admit good isoperimetric constants (\cite{bakry1996levy}), it has been tried to reduce the measure $\mu$ to a more log-concave measure.

Instead of the localization approach, Ronen Eldan~\cite{eldan2013thin} introduced the \textit{stochastic localization} approach.
The idea is to change the measure $\mu$ with a stochastic process over time $t$. This yields a measure decomposition of $\mu$ into a mixture of Gaussians $\mu_t$ with isoperimetric constants $\sqrt{t}$ improving over time. The process is then stopped after a short amount of time $T$, such that $\frac{\mu_T(S)}{\mu(S)}$ is constant with high probability for sets with ${\mu(S)}=\frac{1}{2}$, i.e. it preserves constant mass in large sets. Hence, the isoperimetric constant can be bounded by $\sqrt{T}$ and a proof then breaks down to find the best time $T$ that guarantees the boundedness of $\frac{\mu_T(S)}{\mu(S)}$.

Eldan~\cite{eldan2013thin} proved $\psi_n \geq Cn^{-\frac{1}{3}}\log(n)^{-\frac{1}{2}}$, which was later improved by Lee and Vempala~\cite{lee2017eldan} to $\psi_n \geq Cn^{-\frac{1}{4}}$ and believed to be optimal for some time.

This changed due to a breakthrough, achieved by Yuansi Chen~\cite{chen2021almost} proving
\begin{align*}
    \psi_n \geq C\exp(-\sqrt{\log(n)\log(\log(n))})
\end{align*}
using stochastic localization and a bootstrapping argument, breaking the polynomial bottleneck. This was recently improved by Klartag~\cite{klartag2023logarithmic} using an \textit{improved Lichnerowicz inequality} to obtain 
\begin{align*}
    \psi_n \geq C \log(n)^{-\frac{1}{2}},
\end{align*}
which is currently the best known lower bound.

The KLS Conjecture can also be related to two other conjectures which have been resolved in the past months:
The Thin Shell Conjecture was proven by Klartag and Lehec~\cite{klartag2025thin} and states the following: 
\begin{theorem}
    For any isotropic log-concave measure $\mu$ on $\mathbb{R}^n$, the mass concentrates around a thin spherical shell of radius $\sqrt{n}$ and constant width, i.e. there exists a universal constant $C>0$ such that
\begin{align*}
    \mathbb{E}((\|X\|_2 - \sqrt{n})^2) \leq C .
\end{align*}
\end{theorem}
The fact that the KLS conjecture implies the Thin Shell Conjecture can be readily seen by applying the Poincaré inequality. Remarkably, Eldan~\cite{eldan2013thin} proved the reverse up to $\sqrt{\log(n)}$ using stochastic localization. In particular, the resolution of the Thin Shell Conjecture in 2025 (\cite{klartag2025thin}) also implies the current best bound for the KLS Conjecture.

Another famous problem from high-dimensional convex geometry is Bourgain's slicing problem~\cite{bourgain1986high} from 1986 (now a theorem of Klartag and Lehec~\cite{klartag2025affirmative}):
\begin{theorem}
    There exists a universal constant $C>0$ such that for every dimension $n$ and every convex body $K \subset \mathbb{R}^n$ with unit volume there exists a hyperplane $H \subset \mathbb{R}^n$ such that $\mathrm{Vol}_{n-1}(K \cap H) \geq C$.
\end{theorem}
Moreover, it can be reformulated (\cite{ball1988logarithmically}) for the class of log-concave isotropic measures with density $f$: There is a universal constant $C>0$ such that $f(0)^\frac{1}{n} \leq C$.
Remarkably, the slicing problem has been posed before the KLS Conjecture and the connection was not known for a long time: Eldan and Klartag~\cite{eldan2011approximately} used the log-Laplace transform to show that the Thin Shell Conjecture implies an affirmative answer to Bourgain's problem and therefore the KLS Conjecture also implies Bourgain's slicing conjecture. Based on a new technique due to Guan~\cite{guan2024note}, the problem was first resolved in December 2024 by Klartag and Lehec~\cite{klartag2025affirmative}, before the same authors proved the more general Thin Shell bound a few months later.

\section{Sharp Bounds for Graph Matrices (PNN)}

Given a computation problem on a random input, one natural question is whether there is \textit{some} polynomial time algorithm solving it with high probability. To argue that such a problem is hard, we usually compare it against a set of algorithms that are believed to be strong. Whenever the problem can be phrased as a polynomial optimization problem, a good choice of algorithms is the Sum-of-Squares (SoS) hierarchy. This hierarchy offers an increasing sequence of convex relaxations to the polynomial optimization problem, indexed by the degree of polynomials each relaxation uses. If this degree is $d$, it is known that a solution can be found in $n^{O(d)}$ time using semidefinite programming (SDP). The past decade has seen vast interest in studying this hierarchy as a benchmark, due to the following two reasons.
\begin{enumerate}
    \item For many problems, efficient algorithms were constructed based on it \cite{goemans1995improved, arora2004expander, hopkins2015tensor}, some of which are still state-of-the-art, and for some problems, they have even been shown to be optimal among all SDP relaxations \cite{lee2015lower}.
    \item There has been a well developed machinery to  show lower bounds against SoS hierarchies, meaning that a certain degree $d$ relaxation fails to certify with high probability.
\end{enumerate}

In this post, we focus on the technical aspects of the latter point, for which we highly recommend the survey by Potechin and Rajendran \cite{potechin2023machineryprovingsumofsquareslower}.

Suppose one wants to lower bound the value of the polynomial optimization problem
\begin{align*}
    \min_{x\in S}&\ P(x)\\
    \text{s.t.}&\ g_1(x) = 0, \ldots, g_m(x) = 0
\end{align*}
by some constant $c \in \mathbb{R}$. Instead of checking all instances of $x\in S$, an SoS algorithm tries to find polynomials $f_1(x), \ldots, f_m(x)$ and $h_1(x),\ldots,h_n(x)$ for which the following eponymous identity holds:
\begin{equation*}
    P(x) = c + \sum_{i=1}^m f_i(x)g_i(x) + \sum_{j=1}^n h_j(x)^2.
\end{equation*}
Using degree $d$ relaxation here means that all summands are  polynomials of degree up to $d$. To prevent the SoS algorithm in its quest, one can try finding a linear functional $\tilde{\mathbb{E}}$ on the polynomials of degree at most $d$, called \textit{pseudo-expectation}, satisfying
\begin{enumerate}
    \item (normalization) $\tilde{\mathbb{E}}[1] = 1$,
    \item (feasibility) $\tilde{\mathbb{E}}[f g_i] = 0$ for every $i = 1,\ldots,m$ and polynomial $f$ with $\deg(f g_i) \leq d$,
    \item (positivity) $\tilde{\mathbb{E}}[f^2] \geq 0$ for all polynomials $f$ with $\deg(f) \leq d/2$.
\end{enumerate}
Interestingly, by duality, if there is no such $\tilde{\mathbb{E}}$, the SoS algorithm will succeed. However, finding such $\tilde{\mathbb{E}}$ can be challenging, particularly because the feasibility constraint depends on the random input of the certification problem. A commonly used heuristic is the so-called \textit{pseudo-calibration} \cite{barak2019nearly}, and once it finds some $\tilde{\mathbb{E}}$, the hardest step is checking positivity. We rephrase this condition in terms of the positive semi-definiteness of the \textit{moment matrix} $\Lambda$.

\begin{definition}
    Given a degree $d$ pseudo-expectation $\tilde{\mathbb{E}}$, its associated moment matrix $\Lambda$ has rows and columns indexed by monomials $p$ and $q$ of degree at most $d/2$, and its entries are given by
    \begin{equation*}
        \Lambda[p,q] = \tilde{\mathbb{E}}[pq].
    \end{equation*}
\end{definition}

To deduce $\Lambda \succeq 0$, one needs to understand the spectrum of a random matrix $\Lambda$, whose entries depend on a random input. In particular, the pseudo-calibration heuristic constructs $\Lambda$ whose entries \textit{polynomially} depend on a random input (the degree of such polynomial can be much larger than the SoS degree $d$). We illustrate this statement in the context of the \textit{planted clique problem}.
\begin{center}
    \textit{Given $G \sim G(n,1/2)$ (Erdős–Rényi model), can one certify using degree $d$ SoS that $G$ has no clique of size $k$?}
\end{center}
If we encode the edges of $G$ by i.i.d. Rademachers $(\varepsilon_{i,j})_{i<j}$ ($\pm 1$ for occurrence/absence), then $\Lambda$ (by construction) has entries that are polynomials of those $\binom{n}{2}$ Rademacher variables. At the time, studying the spectrum of such \textit{locally random} matrices was a bottleneck, resolved in \cite{meka2015sum}, and resulted in the development of \textit{graph matrices} \cite{ahn2016graph}. We define them in the restricted setting, using the already-mentioned $(\varepsilon_{i,j})_{i<j}$ as randomness.

\begin{definition}[Graph matrices]\label{def:graphmatrices}
    Let $\alpha$ be a \textit{shape} (a small fixed graph) in which we partition the vertices $V(\alpha)$ into the left ($U_\alpha$) and right ($V_\alpha$) sides, and let $n$ be a large integer. A \emph{realization} is any injective map $\varphi \colon\! V(\alpha) \to [n]$ from the shape vertices to the ground set $[n]$. The \emph{graph matrix} $M_\alpha$ is the $n^{|U_\alpha|} \times n^{|V_\alpha|}$ random matrix, whose rows and columns are indexed by ordered subsets of $[n]$ with cardinalities $|U_\alpha|$ and $|V_\alpha|$, respectively, given by
    \begin{equation*}\label{eq:graphmatrixdef}
        M_\alpha = \sum_{\text{realization } \varphi} \left(\prod_{(i,j) \in E(\alpha)} \varepsilon_{\varphi(i),\varphi(j)}\right) e_{\varphi(U_\alpha)} e_{\varphi(V_\alpha)}^\top.
    \end{equation*}
\end{definition}

Setting the technical aspects of this definition aside, there are two key takeaways.
\begin{enumerate}
    \item Since $\Lambda - \mathbb{E}\Lambda$ is expressible as a linear combination of graph matrices, the spectrum can be understood by examining each graph matrix summand. At a high level, if $\|\Lambda - \mathbb{E}\Lambda\| < s_{min}(\mathbb{E}\Lambda)$, positive semi-definiteness of $\Lambda$ follows; further refinements of this argument across different subspaces of $\Lambda$ yield tighter bounds.
    
    \item Due to the underlying combinatorial symmetry, one can view graph matrices as generalized Wigner models. The classical moment method can be adapted to produce (up to a polylogarithmic factor in $n$) matching norm bounds for any graph matrix $M_\alpha$. Namely, there are functions $f$ and $g$ depending only on the shape $\alpha$, such that
    \begin{equation}\label{eq:gmtx_bound}
        \Omega\left(n^{f(\alpha)}\right) \leq \mathbb{E}\|M_\alpha\| \leq O\left(n^{f(\alpha)}\,(\log n)^{g(\alpha)}\right).
    \end{equation}
\end{enumerate}

Although this recipe has produced many successful results regarding SOS lower bounds over the past decade, the moment method provides only partial insight into the graph matrices themselves. Determining the spectrum \cite{cai2020spectrum} or establishing tighter norm bounds \cite{hsieh2023ellipsoid} has required quite challenging adaptations.

It turns out that graph matrices are only a specific instance of a much broader model, called \textit{matrix chaoses}, that has been studied in operator theory since the 1980s.
\begin{definition}[Matrix Chaos]
    We say a random matrix $X$ is a matrix chaos of order $q$ if there are $m$ i.i.d. random variables $h_1, \ldots, h_m$ such the entries of $X$ are degree $q$ polynomials of those random variables, i.e. there are $m^q$ deterministic matrices $(A_{i_1,\dots,i_q})$ such that
    \begin{equation*}
        X = \sum_{i_1,\ldots,i_q\in[m]} h_{i_1}\cdots h_{i_q} A_{i_1,\dots,i_q}.
    \end{equation*}
\end{definition}

The linear case ($q=1$) has received immense attention in the past decades. Starting from the noncommutative Khintchine (NCK) inequality of Lust-Piquard and Pisier \cite{LPkhintchine91}, a versatile toolkit of \textit{matrix concentration inequalities} has been developed and used on many different questions in numerical algebra, theoretical computer science, statistics, and other fields \cite{tropp2015introduction}. In recent years, the sharpness of such tools has been studied via idealized models from free probability, giving rise to \textit{strong matrix concentration inequalities} \cite{BBvH-Free, BCSvH-Free2, universalityBRvH24}.

The polynomial case ($q>1$) has appeared in many different settings. Nevertheless, researchers have often relied on developing case-by-case tools. However, Haagerup and Pisier \cite{haagerup1993bounded} observed that after decoupling \cite{de2012decoupling}, one can iterate NCK inequality to bound the spectral norm of $X$.

In our recent paper \cite{bandeira2025matrix}, we use this approach to produce a standardized toolkit of matrix chaos inequalities. As graph matrices are matrix chaoses (also observed in \cite{rajendran2023concentration, tulsiani2024simple}), this toolkit is used to rederive norm bounds \eqref{eq:gmtx_bound}. Furthermore, by iterating sharp matrix concentration inequalities, one can remove unnecessary polylogarithmic factors for some instances of $\alpha$, i.e. showing \eqref{eq:gmtx_bound} with $g(\alpha) = 0$.

However, there are still examples of shapes $\alpha$ for which iterated linear inequalities fail to produce sharp bounds. Moreover, there are examples of shapes whose norm bounds do have additional non-constant factors. We conjecture that given a fixed shape $\alpha$, the norm of the associated graph matrix grows as a product of a polynomial and a polylogarithmic factor (depending only on $\alpha$).

\begin{conjecture}[Sharp graph matrix bounds]
    There are functions $f$ and $g$ depending only on the shape $\alpha$, such that
    \begin{equation}\label{eq:sharp_gmtx_bound}
        \mathbb{E}\|M_\alpha\| = \Theta\left(n^{f(\alpha)}\,(\log n)^{g(\alpha)}\right).
    \end{equation}
\end{conjecture}

There is a twofold motivation in determining the correct $f$ and $g$ in \eqref{eq:sharp_gmtx_bound}.
\begin{enumerate}
    \item While current iterated free probability tools fail to certify all shapes that exhibit noncommutative behavior (producing no additional non-constant factors), a correct characterization would hint at which parameters these inequalities should involve and what causes such behavior in matrix chaos. We can view graph matrices as a family of matrix chaoses that possess enough symmetry for an exact characterization to be achievable, while retaining enough degrees of freedom to manifest diverse matrix chaos behaviors.

    \item There have been many applications relying on \eqref{eq:gmtx_bound}, which still have unnecessary polylogarithmic factors in some cases. Any improvement to \eqref{eq:sharp_gmtx_bound} results in better lower bounds for SoS hierarchies. Ultimately, understanding the matrix chaos $\Lambda$ better would provide deeper insight into the pseudo-calibration heuristic.
\end{enumerate}


\bibliographystyle{alpha}
\bibliography{references}

\end{document}